\documentclass[11pt]{article}
\usepackage{amssymb,amsmath}
\usepackage{amsfonts}
\usepackage{graphicx}
\usepackage{tikz}
\usepackage{cancel}
\usepackage{todonotes}
\usepackage{enumitem}
\usetikzlibrary{matrix}
\usetikzlibrary{arrows,shapes}

\numberwithin{equation}{section}
\numberwithin{figure}{section}

\newtheorem{Theorem}{Theorem}[section]
\newtheorem{Prop}[Theorem]{Proposition}

\newtheorem{Lemma}[Theorem]{Lemma}
\newtheorem{Cor}[Theorem]{Corollary}
\newtheorem{Def}[Theorem]{Definition}

\def\tto{\rightrightarrows}

\def\st{\stackrel}

\def\R{\mathbb{R}}

\def\B{\mathbb{B}}
\def\O{\Omega}

\def\dom{\mathrm{dom}\,}
\def\gph{\mathrm{gph}\,}

\def\lm{\lambda}
\def\vt{\vartheta}

\def\al{\alpha}

\def\gg{\gamma}
\def\ph{\varphi}

\def\ox{\bar{x}}

\def\oy{\bar{y}}

\def\dd{\delta}

\def\ve{\varepsilon}
\def\epsilon{\ve}

\def\dom{{\rm dom}\,}
\def\co{{\rm co}\,}

\def\Th{\Theta}
\def\dn{\downarrow}
\def\oR{\overline{\mathbb{R}}}

\def\emp{\emptyset}
\def\epi{{\rm epi}\,}
\def\la{\langle}

\def\ra{\rangle}
\def\tilde{\widetilde}

\def\disp{\displaystyle}

\def\h{\hfill\Box}
\begin{document}
\begin{center}
\textbf{BILEVEL OPTIMIZATION AND VARIATIONAL ANALYSIS}\\[2ex]
BORIS S. MORDUKHOVICH\footnote{Department of Mathematics, Wayne State University, Detroit, Michigan, 48202, USA (boris@math.wayne.edu). This research was partly supported by the USA National Science Foundation under grants DMS-1512846 and DMS-1808978, by the USA Air Force Office of Scientific Research under grant 15RT04, and by Australian Research Council Discovery Project DP-190100555.}\\[2ex]
\end{center}
{\bf Abstract.} This chapter presents a self-contained approach of variational analysis and generalized differentiation to deriving necessary optimality in problems of bilevel optimization with Lipschitzian data. We mainly concentrate on optimistic models, although the developed machinery also applies to pessimistic versions. Some open problems are posed and discussed.\\[1ex]

\section{Introduction}\label{intro}
\setcounter{equation}{0}

Bilevel optimization has been well recognized as a theoretically very challenging and practically important area of applied mathematics. We refer the reader to the monographs \cite{demp03,demp-kalash15,m18}, the extensive bibliographies and commentaries therein, as well as to the advanced material included in this book for various approaches, theoretical and numerical results, and a variety of practical applications of bilevel optimization and related topics.

One of the characteristic features of bilevel optimization problems is their intrinsic {\em nonsmoothness}, even if their initial data are described by linear functions. This makes natural to develop an approach of modern {\em variational analysis and generalized differentiation} to the study and applications of major models in bilevel optimization. It has been done in numerous publications, which are presented and analyzed in the author's recent book \cite{m18}.

The main goal we pursue here is to overview this approach together with the corresponding machinery of variational analysis and to apply it to deriving necessary optimality conditions in optimistic bilevel models with {\em Lipschitzian} data while also commenting on other versions in bilevel optimization with posting open questions. To make this chapter largely {\em self-contained} and more accessible for the reader, we present here the basic background from variational analysis and generalized differentiation, which is needed for applications to bilevel optimization. For brevity and simplicity we confine ourselves to problems in finite-dimensional spaces.

The rest of this work is organized as follows. In Section~\ref{basic} we recall those constructions of {\em generalized differentiation} in variational analysis, which are broadly used in the subsequent text. Section~\ref{extremal} presents the fundamental {\em extremal principle} that is behind generalized differential calculus and applications to optimization in the geometric approach to variational analysis developed in \cite{m06,m18}. Section~\ref{calculus} is devoted to deriving---via the extremal principle---the two basic {\em calculus rules}, which are particularly useful for applications to optimality conditions. In Section~\ref{marginal} we establishing  {\em subdifferential evaluations} and efficient conditions that ensure the local {\em Lipschitz continuity} of {\em optimal value function} in general problems of parametric optimization. These results are crucial for variational applications to bilevel programming.

To proceed with such applications, we first consider in Section~\ref{program} problems of {\em nondifferentiable programming} with Lipschitzian data. Subdifferential necessary optimality conditions for Lipschitzian programs are derived there by using the extremal principle and calculus rules. Section~\ref{sec:bilevel} contains the formulation of the bilevel optimization problems under consideration and the description of the {\em variational approach} to their study. Based on this approach and subdifferentiation of the optimal value functions for lower-level problems, we establish in Section~\ref{optimality} necessary optimality conditions for {\em Lipschitzian bilevel programs}. The other developments in this direction for bilevel optimization problems with Lipschitzian data is presented in Section~\ref{difference} by using the {\em subdifferential difference rule} based on a certain variational technique. The concluding Section~\ref{conclusion} discusses {\em further perspectives} of employing concepts and techniques of variational analysis to bilevel optimization with formulations of some {\em open questions}.

Throughout this chapter we use the standard notation and terminology of variational analysis and generalized differentiation; see, e.g., \cite{m06,m18,rw}.

\section{Basic Constructions of Generalized Differentiation}\label{basic}
\setcounter{equation}{0}

Here we present the basic definitions of generalized normals to sets, coderivatives of set-valued mappings, and subgradients of extended-real-valued functions initiated by the author \cite{m76} that are predominantly used in what follows. The reader is referred to the books \cite{m06,m18,rw} for more details. Developing a {\em geometric approach} to generalized differentiation, we start with normals to sets, then continue with coderivatives of (set-valued and single-valued) mappings, and finally pass to subgradients of extended-real-valued functions.

Given a nonempty set $\O\subset\R^n$, we always suppose without loss of generality that it is locally closed around the reference point $\ox\in\O$. For each $x\in\R^n$ close to $\ox$ consider its (nonempty) Euclidean projector to $\O$ defined by
\begin{equation*}
\Pi(x;\O):=\big\{w\in\O\big|\;\|x-w\|=\disp\min_{u\in\O}\|x-u\|\big\}.
\end{equation*}
Then the (basic, limiting. Mordukhovich) {\em normal cone} to $\O$ at $\ox$ is
\begin{equation}\label{nc}
\begin{array}{ll}
N(\ox;\O):=\Big\{v\in\R^n\Big|&\exists\,x_k\to\ox,\;\exists\,w_k\in\Pi(x_k;\O),\;\exists\,\al_k\ge 0\\
&\mbox{such that }\;\al_k(x_k-w_k)\to v\;\mbox{ as }\;k\to\infty\Big\}.
\end{array}
\end{equation}
The normal cone \eqref{nc} is always closed while may be nonconvex in standard situations; e.g., when $\O$ is the graph of the simplest nonsmooth convex function $|x|$ at $\ox=(0,0)\in\R^2$. Nevertheless, this normal cone and the associated coderivatives of mappings and subdifferentials of functions enjoy comprehensive calculus rules due to variational/extremal principles of variational analysis. Note that $N(\ox)\ne\{0\}$ if and only if $\ox$ is a boundary point of $\O$.

There is a useful representation of the normal cone \eqref{nc} in terms of convex collections of (pre)normal vectors to $\O$ st point nearby $\ox$. Given $x\in\O$ close to $\ox$, we {\em prenormal cone} to $\O$ at $x$ (known also as the regular or Fr\'echet normal cone) is defined by
\begin{equation}\label{pre}
\Hat N(x;\O):=\Big\{v\in\R^n\Big|\;\disp\limsup_{u\st{\O}{\to}x}\frac{\la v,u-x\ra}{\|u-x\|}\le 0\Big\},
\end{equation}
where the symbol $u\st{\O}{\to}x$ means that $u\to x$ with $u\in\O$. Then the prenormal cone \eqref{pre} is always closed and convex while may collapse to $\{0\}$ at boundary points of closed sets, which in fact contradicts the very meaning of generalized normals. If $\O$ is convex, then both normal and prenormal cones reduce to the normal cone of convex analysis. In general we have
\begin{equation}\label{nor-repr}
N(\ox;\O)=\big\{v\in\R^n\big|\;\exists\,x_k\st{\O}{\to}\ox,\;v_k\in\Hat N(x_k;\O)\;\mbox{ with }\;v_k\to v\;\mbox{ as }\;k\to\infty\big\}.
\end{equation}
Note that the limiting representation \eqref{nor-repr} keeps holding if the prenormal cone \eqref{pre} therein is expanded to its {\em $\ve_k$-enlargements} $\Hat{N}_{\ve_k}$ as $\ve_k\dn 0$, where the latter expansions are defined by replacing $0$ with $\ve_k$ on the right-hand side of \eqref{pre}.

Let $F\colon\R^n\tto\R^m$ be a set-valued mapping/multifunction with the values $F(x)\subset\R^m$ and with its graph defined by
\begin{equation*}
\gph F:=\big\{(x,y)\in\R^n\times\R^m\big|\;y\in F(x)\big\}.
\end{equation*}
When $F$ is single-valued, we use use the standard notation $F\colon\R^n\to\R^m$. Assuming that the graph of $F$ is locally closed around $(\ox,\oy)\in\gph F$, we define the {\em coderivative} of $F$ at this point via the normal cone \eqref{nc} to the graph of $F$ by
\begin{equation}\label{cod}
D^*F(\ox,\oy)(w):=\big\{v\in\R^n\big|\;(v,-w)\in N\big((\ox,\oy);\gph F\big)\big\},\quad w\in\R^m.
\end{equation}
Thus $D^*F(\ox,\oy)\colon\R^m\tto\R^n$ is a set-valued positively homogeneous mapping, which reduces to the adjoint/transposed Jacobian for
single-valued mappings $F\colon\R^n\to\R^m$ that are smooth around $\ox$, where $\oy=F(\ox)$ is dropped in this case in the coderivative notation:
\begin{equation*}
D^*F(\ox)(w)=\big\{\nabla F(\ox)^*w\big\}\;\mbox{ for all }\;w\in\R^m.
\end{equation*}
Besides a full calculus available for the coderivative \eqref{cod}, this construction plays an important role in variational analysis and its applications since it provides complete characterizations of fundamental {\em well-posedness} properties of multifunctions concerning Lipschitzian stability, metric regularity, and linear openness/covering. In this work we deal with the {\em Lipschitz-like} (Aubin, pseudo-Lipschitz) property of $F\colon\R^n\tto\R^m$ around $(\ox,\oy)\in\gph F$ defined as follows: there exist neighborhoods $U$ of $\ox$ and $V$ of $\oy$ and a constant $\ell\ge 0$ such that
\begin{equation}\label{lip}
F(x)\cap V\subset F(u)+\ell\|x-u\|\B\;\mbox{ for all }\;x,u\in U,
\end{equation}
where $\B$ stands for the closed unit ball of the space in question. If $V=\R^m$ in \eqref{lip}, then it reduces to the classical local Lipschitzian property of $F$ around $\ox$. The {\em coderivative characterization} of \eqref{lip}, which is called in \cite{rw} the Mordukhovich criterion, tells us that $F$ is Lipschitz-like around $(\ox,\oy)$ if and only if we have
\begin{equation}\label{cc}
D^*F(\ox,\oy)(0)=\{0\}.
\end{equation}
Furthermore, the exact bound (infimum) of all the Lipschitz constant $\{\ell\}$ in \eqref{lip} is calculated as the norm $\|D^*F(\ox,\oy)\|$ of the positively homogeneous coderivative mapping $w\mapsto D^*F(\ox,\oy)(w)$. The reader can find in \cite{m93,m06,m18,rw} different proofs of this result with numerous applications.

Consider finally an extended-real-valued function $\ph\colon\R^n\to\oR:=(-\infty,\infty]$ finite at $\ox$ and lower semicontinuous (l.s.c.) around this point. Denote by
\begin{equation*}
\dom\ph:=\big\{x\in\R^n\big|\;\ph(x)<\infty\big\}\;\mbox{ and }\;\epi\ph:=\big\{(x,\mu)\in\R^n\times\R\big|\;\mu\ge\ph(x)\big\}
\end{equation*}
the domain and epigraph of $\ph$, respectively. Given $\ox\in\dom\ph$ and using the normal cone \eqref{nc} to the epigraph of $\ph$ at $(\ox,\ph(\ox))$, we define the two types of the subdifferentials of $\ph$ at $\ox$: the {\em basic subdifferential} and the {\em singular subdifferential} by, respectively,
\begin{equation}\label{bs}
\partial\ph(\ox):=\big\{v\in\R^n\big|\;(v,-1)\in N\big((\ox,\ph(\ox));\epi\ph\big)\big\},
\end{equation}
\begin{equation}\label{ss}
\partial^\infty\ph(\ox):=\big\{v\in\R^n\big|\;(v,0)\in N\big((\ox,\ph(\ox));\epi\ph\big)\big\}.
\end{equation}
The basic subdifferential \eqref{bs} reduces to the gradient $\{\nabla\ph(\ox)\}$ for smooth functions and to the subdifferential of convex analysis if $\ph$ is convex. Observe that $\partial\ph(\ox)=D^*E_\ph(\ox,\ph(\ox))(1)$ and $\partial^\infty\ph(\ox)=D^*E_\ph(\ox,\ph(\ox))(0)$ via the coderivative \eqref{cod} of the epigraphical multifunction $E_\ph\colon\R^n\tto\R$ defined by $E_\ph(x):=\{\mu\in\R|\;\mu\ge\ph(x)\}$. Thus the coderivative characterization \eqref{cc} of the Lipschitz-like property of multifunctions implies that a lower semicontinuous function $\ph$ is {\em locally Lipschitzian} around $\ox$ if and only if
\begin{equation}\label{lip-func}
\partial^\infty\ph(\ox)=\{0\}.
\end{equation}
Note also that, given any (closed) set $\O\subset\R^n$ with its {\em indicator function} $\dd(x;O)=\dd_\O(x)$ equal $0$ for $x\in\O$ and $\infty$ otherwise, we have that
\begin{equation}\label{indic}
\partial\dd(\ox;\O)=\partial^\infty\dd(\ox;\O)=N(\ox;\O)\;\mbox{ whenever }\;\ox\in\O.
\end{equation}
Both subdifferentials \eqref{bs} and \eqref{ss} admit limiting representations in terms of the {\em presubdifferential}, or {\em regular subdifferential}
\begin{equation}\label{rs}
\Hat\partial\ph(x):=\Big\{v\in\R^n\Big|\;\disp\frac{\ph(u)-\ph(x)-\la v,u-x\ra}{\|u-x\|}\ge 0\Big\}
\end{equation}
of $\ph$ at points $x$ close to $\ox$. Namely, we have
\begin{equation}\label{bs1}
\partial\ph(\ox)=\big\{v\in\R^n\big|\;\exists\,x_k\st{\ph}{\to}\ox,\;\exists\,v_k\to v\;\mbox{ with }\;v_k\in\Hat\partial\ph(x_k)\;\mbox{ as }\;k\to\infty\big\},
\end{equation}
\begin{equation}\label{ss1}
\partial^\infty\ph(\ox)=\big\{v\in\R^n\big|\;\exists\,x_k\st{\ph}{\to}\ox,\;\exists\lm_k\dn 0,\;\exists\,v_k\to v\;\mbox{ with }\;v_k\in\lm_k\Hat\partial\ph(x_k)\;\mbox{ as }\;k\to\infty\big\},
\end{equation}
where the symbol $x\st{\ph}{\to}\ox$ indicates that $x\to\ox$ with $\ph(x)\to\ph(\ox)$. Note that the presubdifferential \eqref{rs} is related to the prenormal cone \eqref{pre} as in \eqref{bs} and is also used in variational analysis under the names of the Fr\'echet subdifferential and the viscosity subdifferential.
Similarly to the case of basic normals in \eqref{nor-repr} it is not hard to observe that we still have the subdifferential representations in \eqref{bs1} and \eqref{ss1} if the presubdifferential \eqref{rs} therein is expanded by its {\em $\ve_k$-enlargements} $\Hat\partial_{\ve_k}\ph$ defined with replacing $0$ on the right-hand side of \eqref{rs} by $-\ve_k$.

\section{Extremal Principle in Variational Analysis}\label{extremal}
\setcounter{equation}{0}

In this section we recall, following \cite{km}, the notion of locally extremal points for systems of finitely many sets and then derive the fundamental {\em extremal principle}, which gives us necessary conditions for extremality of closed set systems in $\R^n$.

\begin{Def}\label{set-ext} Let $\O_1,\ldots,\O_s$ as $s\ge 2$ be nonempty subsets of $\R^n$, which are assumed to be locally closed around their common point $\ox$. We say that $\ox$ is a {\sc locally extremal point} of the set system $\{\O_1,\ldots,\O_s\}$ if there exist a neighborhood $U$ of $\ox$ and sequences of vectors $a_{ik}\in\R^n$, $i=1,\ldots,s$, such that $a_{ik}\to 0$ as $k\to\infty$ for all $i\in\{1,\ldots,s\}$ and
\begin{equation}\label{ext1}
\bigcap_{i=1}^s\big(\O_i-a_{ik}\big)\cap U=\emp\;\mbox{ whenever }\;k=1,2,\ldots.
\end{equation}
\end{Def}
Observe that for the case of two sets $\O_1,\O_2$ containing $\ox$ the above definition can be equivalently reformulated as follows: there is a neighborhood $U$ of $\ox$ such that for any $\ve>0$ there exists a vector $a\in\R^n$ with $\|a\|\le\ve$ and $(\O_1-a)\cap\O_2\cap U=\emp$.

It is easy to see of a closed set $\O$ an its boundary point $\ox$ form the extremal system $\{\O,\{\ox\}\}$. Furthermore, the introduced notion of set extremality covers various notions of optimality and equilibria in problems of scalar and vector optimization. In particular, a local minimizer $\ox$ of the general constrained optimization problem
\begin{equation*}
\mbox{minimize }\;\ph(x)\;\mbox{ subject to }\;x\in\O\subset\R^n,
\end{equation*}
where $\ph$ is l.s.c.\ and $\O$ is closed around $\ox$, corresponds to the locally extremal point $(\ox,\ph(\ox))$ of the sets $\O_1:=\epi\ph$ and $\O_2:=\O\times\{\ph(\ox)\}$. As we see below, extremal systems naturally arise in deriving calculus rules of generalized differentiation.

Now we are ready to formulate and prove the basic extremal principle of variational analysis for systems of finitely many closed sets in $\R^n$ by using the normal cone construction \eqref{nc}.

\begin{Theorem}\label{ep} Let $\ox$ be a locally extremal point of the system $\{\O_1,\ldots,\O_s\}$ of nonempty subsets of $\R^n$, which are locally closed around $\ox$. Then there exist generalized normals $v_i\in N(\ox;\O_i)$ for $i=1,\ldots,s$, not equal to zero simultaneously, such that we have the generalized Euler equation
\begin{equation}\label{euler}
v_1+\ldots+v_s=0.
\end{equation}
\end{Theorem}
{\bf Proof.} Using Definition~\ref{set-ext}, suppose without loss of generality that $U=\R^n$. Taking the sequences $\{a_{ik}\}$ therein, for each $k=1,2,\ldots$ consider the unconstrained optimization problem:
\begin{equation}\label{dk}
\mbox{minimize }\;\ph_k(x):=\Big[\sum_{i=1}^s d^2(x+a_{ik};\O_i)\Big]^{1/2}+\|x-\ox\|^2,\quad x\in\R^n,
\end{equation}
where $d(x;\O)$ indicates the Euclidean distance between $x$ and $\O$. Since $\ph_k$ is continuous and the level sets of it are bounded, we deduce from the classical Weierstrass theorem that there exists an optimal solution $x_k$ to each problem \eqref{dk} as $k=1,2,\ldots$. It follows from the crucial extremality requirement \eqref{ext1} in Definition~\ref{set-ext} that
\begin{equation}\label{gammak}
\gg_k:=\Big[\sum_{i=1}^sd^2(x_k+a_{ik};\O_i)\Big]^{1/2}>0.
\end{equation}
The optimality of $x_k$ in \eqref{dk} tells us that
\begin{equation*}
\ph_k(x_k)=\gg_k+\|x_k-\ox\|^2\le\Big[\sum_{i=1}^s\|a_{ik}\|^2\Big]^{1/2}\dn 0,
\end{equation*}
and so $\gg_k\dn 0$ and $x_k\to\ox$ as $k\to\infty$. By the closedness of the sets $\O_i$, $i=1,\ldots,s$, around $\ox$, we pick $w_{ik}\in\Pi(x_k+a_{ik};\O_i$ and for each $k$ form another unconstrained optimization problem:
\begin{equation}\label{ext2}
\mbox{minimize }\;\psi_k(x):=\Big[\sum_{i=1}^s\|x+a_{ik}-w_{ik}\|^2\Big]^{1/2}+\|x-\ox\|^2,\quad x\in\R^n
\end{equation}
which obviously has the same optimal solution $x_k$. In contrast to $\ph_k$ in \eqref{dk}, the function $\psi_k$ in \eqref{ext2} is {\em differentiable} at $x_k$ due to \eqref{gammak}, Thus applying the Fermat rule in \eqref{ext2} tells us that
\begin{equation}\label{ext3}
\nabla\psi_k(x_k)=\sum_{i=1}^s v_{ik}+2(x_k-\ox)=0
\end{equation}
with $v_{ik}:=(x_k+a_{ik}-w_{ik})/\gg_k$, $i=1,\ldots,s$, satisfying
\begin{equation}\label{ext4}
\|v_{1k}\|^2+\ldots+\|v_{sk}\|^2=1\;\mbox{ for all }\;k=1,2,\ldots.
\end{equation}
Remembering the compactness of the unit sphere in $\R^n$, we get by passing to the limit as $k\to\infty$ in \eqref{ext3} and \eqref{ext4} that there exist $v_1,\ldots,v_s$, not equal to zero simultaneously, for which \eqref{euler} holds. Finally, it follows directly from the above constructions and the normal cone definition \eqref{nc} that $v_i\in N(\ox;\O_i)$ for all $i=1,\ldots,s$. This completes the proof of the theorem. $\h$\vspace*{0.05in}

Since for convex sets $\O$ the normal cone \eqref{nc} reduces to the normal cone of convex analysis
\begin{equation*}
N(\ox;\O):=\big\{v\in\R^n\big|\;\la v,x-\ox\ra\le 0\;\mbox{ whenever }\;x\in\O\big\},
\end{equation*}
the extremal principle of Theorem~\ref{ep} can be treated as a {\em variational extension} of the classical separation theorem to the case of finitely many {\em nonconvex} sets in $\R^n$.

\section{Fundamental Calculus Rules}\label{calculus}
\setcounter{equation}{0}

Employing the extremal principe, we derive here two fundamental rules of generalized differential calculus, which are broadly used in this chapter and from which many other calculus rules follow; see \cite{m06,m18}. The first result is the {\em intersection rule} for basic normals \eqref{nc}.

\begin{Theorem}\label{nor-inter} Let $\O_1,\ldots,\O_s$ be nonempty subsets of $\R^n$, which are locally closed around their common point $\ox$. Assume the validity of the following qualification condition:
\begin{equation}\label{qc}
\big[x_i\in N(\ox;\O_i),\;x_1+\ldots+x_s=0\big]\Longrightarrow x_i=0\;\mbox{ for all }\;i=1,\ldots,s.
\end{equation}
Then we have the normal cone intersection rule
\begin{equation}\label{inter-rule1}
\disp N\Big(\ox;\bigcap_{i=1}^s\O_i\Big)\subset N(\ox;\O_1)+\ldots+N(\ox;\O_s).
\end{equation}
\end{Theorem}
{\bf Proof.} Arguing by induction, we first verify the result for $s=2$, Pick any $v\in N(\ox;\O_1\cap\O_2)$ and use the normal cone representation \eqref{nor-repr}. It gives us sequences $x_k\to\ox$ with $x_k\in\O_1\cap\O_{2k}$ and $v_k\to v$ with $v_k\in\Hat N(\ox,\O_1\cap\O_2)$ as $k\to\infty$. Take any sequence $\ve_k\dn 0$ and construct the sets
\begin{eqnarray*}
\Th_1:=\O_1\times\R_+,\quad\Th_{2k}:=\big\{(x,\al)\big|\;x\in\O_2,\,\la v_k,x-x_k\ra-\ve_k\|x-x_k\|\ge\al\big\}\;\mbox{ for any }\;k=1,2,\ldots,
\end{eqnarray*}
These sets are obviously closed around $(x_k,0)\in\Th_1\cap\Th_{2k}$ for all $k$ sufficiently large. Furthermore, it follows from the prenormal cone definition \eqref{pre} that there exists a neighborhood $U$ of $x_k$ with
\begin{equation*}
\Th_1\cap\big(\Th_{2k}-(0,\gg)\big)\cap(U\times\R)=\emp
\end{equation*}
for small numbers $\gg>0$. Thus the pair $(x_k,0)$ is a locally extremal point of the set system $\{\Th_1,\Th_{2k}\}$ for such $k$. Applying to this system the extremal principle from Theorem~\ref{ep} gives us pairs $(u_k,\lm_k)$ from the unit sphere in $\R^{n+1}$ for which
\begin{equation}\label{inter2}
(u_k,\lm_k)\in N\big((x_k,0);\Th_1\big)\;\mbox{ and }\;(-u_k,-\lm_k)\in N\big((x_k,0);\Th_{2k}\big).
\end{equation}
Passing to a subsequence if needed, we get $((u_k,\lm_k)\to(u,\lm)$ as $k\to\infty$ for some $(u,\lm)\in\R^{n+1}$ with $\|(u,\lm)\|=1$. Passing to the limit in the first inclusion of \eqref{inter2} gives us $(u,\lm)\in N((\ox,0);\Th_1)$, which immediately implies that $u\in\O_1$ and $\lm\le 0$. On the other hand, the limiting procedure in the second inclusion of \eqref{inter2} leads us by the structure of $\Th_{2k}$ to
\begin{equation*}
(-\lm v-u,\lm)\in N\big((\ox,0);\O_2\times\R_+\big).
\end{equation*}\
Assuming there that $\lm=0$ contradicts the qualification condition \eqref{qc} for $s=2$. Thus $\lm<0$, which readily implies that $v\in N(\ox;\O_1)+N(\ox;\O_2)$.

To proceed finally by induction for $s>2$, we observe that the induction assumption for \eqref{inter-rule1} in the previous step yields the validity of the qualification condition \eqref{qc} needed for the current step of induction. This completes the proof of the theorem. $\h$\vspace*{0.05in}

Next we derive the {\em subdifferential sum rules} concerning both basic subdifferential \eqref{bs} and singular subdifferential \eqref{ss}. For our subsequent applications to bilevel optimization, it is sufficient to consider the case where {\em all but one} of the functions involved in summation are locally {\em Lipschitzian} around the reference point. This case allows us to obtain the subdifferential sum rules {\em without any} qualification conditions.

\begin{Theorem}\label{sub-sum} Let $\ph_1\colon\R^n\to\oR$ be l.s.c.\ around $\ox\in\dom\ph_1$, and let $\ph_i\colon\R^n\to\oR$ for $i=2,\ldots,s$ and $s\ge 2$ be locally Lipschitzian around $\ox$. Then we have the sum rules
\begin{equation}\label{bs-sum}
\partial\Big(\disp\sum_{i=1}^s\ph_i\Big)(\ox)\subset\sum_{i=1}^s\partial\ph_i(\ox),
\end{equation}
\begin{equation}\label{ss-sum}
\partial^\infty\Big(\disp\sum_{i=1}^s\ph_i\Big)(\ox)=\partial^\infty\ph_1(\ox).
\end{equation}
\end{Theorem}
{\bf Proof.} We consider the case where only two functions are under summation since the general case of finitely many functions obviously follows by induction. Let us start with the basic subdifferential sum rule \eqref{bs-sum} for $s=2$ therein.

Pick any $v\in\partial(\ph_1+\ph_2)(\ox)$ and get by definition \eqref{bs} that
\begin{eqnarray*}
(v,-1)\in N\big((\ox,(\ph_1+\ph_2)(\ox));\epi(\ph_1+\ph_2)\big).
\end{eqnarray*}
Then construct the sets
\begin{equation*}
\O_i:=\big\{(x,\mu_1,\mu_2)\in\R^n\times\R\times\R\big|\;\mu_i\ge\ph_i(x)\big\}\;\mbox{ for }\;i=1,2.
\end{equation*}
Denoting $\bar\mu_i:=\ph_i(\ox)$, $i=1,2$, we obviously have that the sets $\O_1$ and $\O_2$ are locally closed around the triple $(\ox,\bar\mu_1,\bar\mu_2)\in\O_1\cap\O_2$. It is easy to check that  $(v,-1,-1)\in N((\ox,\bar\mu_1,\bar\mu_2);\O_1\cap\O_2)$. Applying now to this set intersection the normal cone intersection rule from Theorem~\ref{nor-inter}, we observe that the qualification condition \eqref{qc} is automatically satisfied in this case due to the singular subdifferential characterization \eqref{lip-func} of the local Lipschitz continuity. Hence we get pairs $(v_i,-\lm_i)\in N((\ox,\bar\mu_i);\epi\ph_i)$ for $i=1,2$ satisfying the condition
\begin{equation*}
(v,-1,-1)=(v_1,-\lm_1,0)+(v_2,0,-\lm_2),
\end{equation*}
which implies that $v=v_1+v_2$ and $\lm_1=\lm_2=-1$. Therefore it shows that $v_i\in\partial\ph_i(\ox)$ for $i=1,2$, and thus the sum rule \eqref{bs-sum} is verified.

Next we proceed with the proof of \eqref{ss-sum} for $s=2$ starting with verifying the inclusion ``$\subset$" therein. Pick $v\in\partial^\infty(\ph_1+\ph_2)(\ox)$ and find by definition sequences $\gg_k\dn 0$, $(x_k,\mu_k)\st{\rm\small
epi(\ph_1+\ph_2)}{\longrightarrow}(\ox,(\ph_1+\ph_2)(\ox))$, $v_k\to v$, $\nu_k\to 0$, and $\eta_k\dn 0$ such that
$$
\la v_k,x-x_k\ra+\nu_k(\mu-\mu_k)\le\gg_k(\|x-x_k\|+|\mu-\mu_k|)
$$
whenever $(x,\mu)\in\epi(\ph_1+\ph_2)$ with $x\in x_k+\eta_k\B$ and $|\mu-\mu_k|\le\eta_k$ as $k=1,2,\ldots$. Taking a Lipschitz constant $\ell>0$ of $\ph_2$ around $\ox$, denote $\Tilde{\eta}_k:=\eta_k/2(\ell+1)$ and $\tilde{\mu}_k:=\mu_k-\ph_2(x_k)$. Then $(x_k,\tilde{\mu}_k)\st{\rm\small epi\ph_1}{\longrightarrow}(\ox,\ph_1(\ox))$ and
$$
(x,\mu+\ph_2(x))\in\epi(\ph_1+\ph_2),\quad|(\mu+\ph_2(x))-\mu_k|\le\eta_k
$$
for all $(x,\mu)\in\epi\ph_1$, $x\in x_k+\tilde{\eta}_k\B$, and
$|\mu-\Tilde{\mu}_k|\le\tilde{\eta}_k$. Therefore
$$
\la v_k,x-x_k\ra+\nu_k(\mu-\Tilde{\mu}_k)\le\ve_k(\|x-x_k\|+|\mu-\Tilde{\mu}_k|)
\mbox{ with }\;\ve_k:=\gg_k(1+\ell)+|\nu_k|\ell
$$
if $(x,\mu)\in\epi\ph_1$ with $x\in x_k+\Tilde{\eta}_k\B$ and $|\mu-\Tilde{\mu}_k|\le\Tilde{\eta}_k$. It yields $(v_k,\nu_k)\in\Hat{N}_{\ve_k}((x_k,\Tilde{\mu}_k);\epi\ph_1)$ for all $k=1,2,\ldots$, and so $(v,0)\in N((\ox,\ph(\ox));\epi\ph_1)$ since $\ve_k\dn 0$ as $k\to\infty$. This verifies the inclusion ``$\subset$" in \eqref{ss-sum}. Applying it to the sum
$\ph_1=(\ph_1+\ph_2)+(-\ph_2)$ yields $\partial^\infty\ph_1(\ox)\subset\partial^\infty(\ph_2+\ph_1)(\ox)$, which justifies the equality in \eqref{ss-sum} and thus completes the proof. $\h$

\section{Subdifferentials and Lipschitz Continuity of Value Functions}\label{marginal}
\setcounter{equation}{0}

In this section we consider the class of extended-real-valued functions $\vt\colon\R^n\to\oR$ defined by
\begin{equation}\label{marg-def}
\vt(x):=\inf\big\{\ph(x,y)\big|\;y\in F(x)\big\},\quad x\in\R^n,
\end{equation}
where $\ph\colon\R^n\times\R^m\to\oR$ is an l.s.c.\ function, and where $F\colon\R^n\tto\R^m$ is a set-valued mapping of closed graph. We can view \eqref{marg-def} as the {\em optimal value function} in the problem of parametric optimization described as follow:
$$
\mbox{minimize }\;\ph(x,y)\;\mbox{ subject to }\;y\in F(x)
$$
with the cost function $\ph$ and the constraint mapping $F$, where $y$ and $x$ are the decision and parameter variables, respectively. Functions of this type are also known in variational analysis under the name of ``marginal functions." A characteristic feature of such functions is their nonsmoothness regardless of the smoothness of the cost function $\ph$ and the simplicity of the constraint mapping $F$ that may nicely behave on the parameter $x$.

As seen below, functions of type \eqref{marg-def} play a crucial role in applications to bilevel optimization while revealing {\em intrinsic nonsmoothness} of the latter class of optimization problems. This section presents evaluations of both basic and singular subdifferentials of \eqref{marg-def}, which are equally important for the aforementioned applications. {\em Singular subdifferential} evaluations are used for establishing the {\em Lipschitz continuity} of $\vt(x)$ with respect to the parameter $x$ that allows us to reduce the bilevel model under consideration to a single-level problem of Lipschitzian programming. On the other hand, {\em basic subdifferential} evaluations open the gate to derive in this way {\em necessary optimality conditions} for Lipschitzian bilevel programs.

To proceed, let us consider the {\em argminimum mapping} $M\colon\R^n\tto\R^m$ associated with by
\begin{equation}\label{arg}
M(x):=\big\{y\in F(x)\big|\;\ph(x,y)=\vt(x)\big\},\quad x\in\R^n,
\end{equation}
and recall that this mapping is {\em inner semicontinuous} at $(\ox,\oy)\in\gph M$ if for every sequence $x_k\st{{\rm\small dom}\,M}{\longrightarrow}\ox$ there exists a sequence $y_k\in M(x_k)$ that converges to $\oy$ as $k\to\infty$. Observe that the inner semicontinuity of $M$ at $(\ox,\oy)$ is implied by its Lipschitz-like property at this point.\vspace*{0.05in}

The following theorem gives us efficient upper estimates of both basic and singular subdifferentials of the optimal value function $\vt$ needed for subsequent applications. We confine ourselves to the case of local Lipschitz continuity of the cost function $\ph$ in \eqref{marg-def} that is sufficient to apply to deriving necessary optimality conditions for bilevel programs in Sections~\ref{optimality} and \ref{difference}.

\begin{Theorem}\label{marg-sub} Let the argminimum mapping \eqref{arg} be inner semicontinuous at $(\ox,\oy)\in\gph M$, and let the cost function $\ph$ be locally Lipschitzian around this point. Then we have
\begin{eqnarray}\label{marg-bs}
\partial\vt(\ox)\subset\bigcup_{(v,w)\in\partial\ph(\ox,\oy)}\Big[v+D^*F(\ox,\oy)(w)\Big],
\end{eqnarray}
\begin{eqnarray}\label{marg-ss}
\partial^\infty\vt(\ox)\subset D^*F(\ox,\oy)(0).
\end{eqnarray}
\end{Theorem}
{\bf Proof.} To start with the verification of \eqref{marg-bs}, consider the extended-real-valued  function $\psi\colon\R^n\times\R^m\to\oR$ defined via the indicator function of the set $\gph F$ by
\begin{equation}\label{psi}
\psi(x,y):=\ph(x,y)+\delta\big((x,y);\gph F\big)\;\mbox{ for all }\;(x,y)\in\R^n\times\R^m
\end{equation}
and prove first the fulfillment of the estimate
\begin{eqnarray}\label{bs-est}
\partial\vt(\ox)\subset\big\{v\in\R^n\big|\;(v,0)\in\partial\psi(\ox,\oy)\big\}.
\end{eqnarray}
Indeed, pick any subgradient $v\in\partial\vt(\ox)$ and get from its representation in \eqref{bs1} sequences $x_k\st{\vt}{\to}\ox$ and $v_k\to v$ with $v_k\in\Hat{\partial}\vt(x_k)$ as $k\to\infty$. Based on definition \eqref{rs}, for any sequence $\ve_k\dn 0$ there exists $\eta_k\dn 0$ as $k\to\infty$ such that $$
\la v_k,x-x_k\ra\le\vt(x)-\vt(x_k)+\ve_k\|x-x_k\|\;\mbox{ whenever }\;x\in x_k+\eta_k\B,\quad k=1,2,\ldots.
$$
This ensures by using the constructions above that
$$
\la(v_k,0),(x,y)-(x_k,y_k)\ra\le\psi(x,y)-\psi(x_k,y_k)+\ve_k\big(\|x-x_k\|+\|y-y_k\|\big)
$$
for all $y_k\in M(x_k)$ and $(x,y)\in(x_k,y_k)+\eta_k\B$. This tells us that $(v_k,0)\in\Hat{\partial}_{\ve_k}\psi(x_k,y_k)$ for all $k=1,2,\ldots$. Employing further the inner semicontinuous of the argminimum mapping $M$ at $(\ox,\oy)$, we find a sequence of $y_k\in M(x_k)$ converging to $\oy$ as $k\to\infty$. It follows from imposed convergence $\vt(x_k)\to\vt(\ox)$ that $\psi(x_k,y_k)\to\psi(\ox,\oy)$. Hence we arrive at $(v,0)\in\partial\psi(\ox,\oy)$ by passing to the limit as $k\to\infty$, which verifies therefore the validity of the upper estimate \eqref{bs-est}. To derive from \eqref{bs-est} the one in \eqref{marg-bs} claimed in the theorem, it remains to use in \eqref{bs-est} the basic subdifferential sum rule \eqref{bs-sum} from Theorem~\ref{sub-sum} combining it with subdifferentiation of the indicator function in \eqref{indic} and the coderivative definition in \eqref{cod}.

Next we verify the singular subdifferential estimate
\begin{equation}\label{ss-est}
\partial^\infty\vt(\ox)\subset\big\{v\in\R^n\big|\;(v,0)\in\partial^\infty\psi(\ox,\oy)\big\}
\end{equation}
for optimal value function \eqref{marg-def} in terms of the auxiliary function \eqref{psi} under the assumptions made. Picking $v\in\partial^\infty\vt(\ox)$ and taking any sequence $\ve_k\dn 0$, find by \eqref{ss1} sequences $x_k\st{\vt}{\to}\ox$, $(v_k,\nu_k)\to(v,0)$, and $\eta_k\dn 0$ as $k\to\infty$ satisfying
$$
\la v_k,x-x_k\ra+\nu_k(\mu-\mu_k)\le\ve_k\big(\|x-x_k\|+|\mu-\mu_k|\big)
$$
for all $(x,\mu)\in\epi\vt$, $x\in x_k+\eta_k\B$, and $|\mu-\mu_k|\le\eta_k$. The assumed inner semicontinuity of \eqref{arg} ensures the existence of sequences $y_k\st{M(x_k)}{\longrightarrow}\oy$ and $\mu_k\dn\psi(\ox)$ such that
$$
(v_k,0,\nu_k)\in\Hat{N}_{\ve_k}\big((x_k,y_k,\mu_k);\epi\psi\big)\;\mbox{ for all }\;k=1,2,\ldots,
$$
via the $\ve_k$-enlargements $\Hat N_{\ve_k}$ of the prenormal cone to the epigraph of $\psi$. This gives us \eqref{ss-est} by passing to the limit as $k\to\infty$. Applying finally to $\partial^\infty\psi$ in \eqref{ss-est} the singular subdifferential relation \eqref{ss-sum} from Theorem~\ref{sub-sum} with taking into account the singular subdifferential calculation in \eqref{indic} together with the coderivative definition \eqref{cod}, we arrive at the claimed upper estimate \eqref{marg-ss} and thus complete the proof of the theorem. $\h$\vspace*{0.05in}

As mentioned above, in our applications to bilevel programming we need to have verifiable conditions that ensure the local Lipschitz continuity of the optimal value function \eqref{marg-def}. This is provided by the following corollary, which is a direct consequence of Theorem~\ref{marg-sub} and the coderivative criterion \eqref{cc} for the Lipschitz-like property.

\begin{Cor}\label{marg-lip1} In addition to the assumptions of Theorem~{\rm\ref{marg-sub}}, suppose that the constraint mapping $F$ is Lipschitz-like around $(\ox,\oy)\in\gph M$ in \eqref{arg}. Then the optimal value function \eqref{marg-def} is locally Lipschitzian around $(\ox,\oy)$.
\end{Cor}
{\bf Proof}. We know from the coderivative criterion \eqref{cc} that $F$ is Lipschitz-like around $(\ox,\oy)$ if and only if $D^*F(\ox,\oy)(0)=\{0\}$. Applying it to \eqref{marg-ss} tells us that the assumed Lipschitz-like property of the constraint mapping $F$ in \eqref{marg-def} ensures that $\partial^\infty\vt(\ox)=\{0\}$. Furthermore, it easily follows from the assumptions made that the optimal value function is l.s,c.\ around $\ox$. Thus $\vt$ is locally Lipschitzian around $\ox$ by the characterization of this property given in \eqref{lip-func}. $\h$

\section{Problems of Lipschitzian Programming}\label{program}
\setcounter{equation}{0}

Before deriving necessary optimality conditions in Lipschitzian problems of bilevel optimization in the subsequent sections, we devote this section to problems of single-level Lipschitzian programming. The results obtained here are based on the extremal principle and subdifferential characterization of local Lipschitzian functions while being instrumental for applications to bilevel programs given in Section~\ref{optimality}.\vspace*{0.05in}

The mathematical program under consideration here is as follows:
\begin{equation}\label{mp}
\begin{array}{ll}
\mbox{minimize }\;\ph_0(x)\;\mbox{ subject to }\\
\ph_i(x)\le 0\;\mbox{ for all }\;i=1,\ldots,m,
\end{array}
\end{equation}
where the functions $\ph_i\colon\R^n\to\R$, $i=0,\ldots,m$, are locally Lipschitzian around the reference point $\ox$. The next theorem provides necessary optimality condition in problem \eqref{mp} of Lipschitzian programming that are expressed in terms if the basic subdifferential \eqref{bs}.

\begin{Theorem}\label{lip-mp} Let $\ox$ be a feasible solution to problem \eqref{mp} that gives a local minimum to the cost function $\ph_0$ therein. then there exist multipliers $\lm_0,\ldots,\lm_m$ satisfying the sign conditions
\begin{equation}\label{sign}
\lm_i\ge 0\;\mbox{ for all }\;i=0,\ldots,m,
\end{equation}
the nontriviality conditions
\begin{equation}\label{nontr}
\lm_0+\ldots+\lm_m\ne 0,
\end{equation}
the complementary slackness conditions
\begin{equation}\label{comp}
\lm_i\ph_i(\ox)=0\;\mbox{ whenever }\;i=1,\ldots,m,
\end{equation}
and the subdifferential Lagrangian inclusion
\begin{equation}\label{lagr}
0\in\disp\sum_{i=0}^{m}\lm_i\partial\ph_i(\ox).
\end{equation}
Assume in addition that
\begin{equation}\label{mfcq}
\Big[\disp\sum_{i\in I(\ox)}\lm_i v_i=0,\;\lm_i\ge 0\Big]\Longrightarrow\Big[\lm_i=0\;\mbox{ for all }\;i\in I(\ox)\Big]
\end{equation}
whenever $v_i\in\partial\ph_i(\ox)$ with $I(\ox):=\big\{i\in\{1,\ldots,m\}\big|\;\ph_i(\ox)=0\big\}$. Then the necessary optimality conditions formulated above hold with $\lm_0=1$.
\end{Theorem}
{\bf Proof.} Supposing without loss of generality that $\ph_0(\ox)=0$, consider the point $(\ox,0)\in\R^n\times\R^m$ and form the following system of $m+1$ sets in the space $\R^n\times\R^m$:
\begin{equation}\label{sets}
\O_i:=\big\{(x,\mu_0,\ldots,\mu_m)\in\R^n\times\R^m\big|\;(x,\mu_i)\in\epi\ph_i\big\}\;\mbox{ for }\;i=0.\ldots,m.
\end{equation}
It ia obvious that $(\ox,0)\in\O_0\cap\ldots\cap\O_m$ and that all the sets $\O_i$, $i=0,\ldots,m$, are locally closed around $(\ox,0)$. Furthermore, there exists a neighborhood $U$ of the local minimizer $\ox$ such that for any $\ve>$ we find $\nu\in(0,\ve)$ ensuring that
\begin{equation}\label{ext2}
\big(\O_1-a\big)\disp\bigcap_{i=1}^m\O_i\cap\big(U\times\{0\}\big)=\emp,
\end{equation}
where $a:=(0,\nu,0,\ldots,0)\in\R^n\times\R^m$ with $\nu\in\R$ standing at the first position after $0\in\R^n$. Indeed, the negation of \eqref{ext2} contradicts the local minimality of $\ox$ in \eqref{mp}. Having \eqref{ext2} gives us \eqref{ext1} for the set system \eqref{sets} and thus verifies that $(\ox,0)$ is a locally extremal point of these sets. Applying now the extremal principle from Theorem~\ref{ep} to $\{\O_0,\ldots,\O_m)$ at $(\ox,0)$ with taking into account the structures of $\O_i$, we get pairs $(v_0,\lm_0),\ldots,(v_m,\lm_m)\in\R^n\times\R$ such that
\begin{equation}\label{ext3}
(v_i,-\lm_i)\in N\big((\ox,0);\epi\ph_i\big)\;\mbox{ for all }\;i=0.\ldots,m,
\end{equation}
\begin{equation}\label{nontr1}
\disp\sum_{i=0}^m\|(v_i,\lm_i)\|\ne 0,
\end{equation}
\begin{equation}\label{eul}
(v_0,-\lm_0)+\ldots+(v_m,-\lm_m)=(0,0),
\end{equation}
It easily follows from \eqref{ext3} and the structure of the epigraphical sets in \eqref{ext3} that the sign conditions \eqref{sign} are satisfied. Furthermore, the singular subdifferential criterion \eqref{lip-func} for the local Lipschitz continuity of the functions $\ph_i$ as $i=0,\ldots,m$ being combined with the sign conditions \eqref{sign} and the definitions \eqref{bs} and \eqref{ss} of the basic and singular subdifferentials, respectively, tells us that the inclusions in \eqref{ext3} are equivalent to
\begin{equation}\label{ext4}
v_i\in\lm_i\partial\ph_i(\ox)\;\mbox{ for all }\;i=0,\ldots,m.
\end{equation}
This ensures that the nontriviality conditions in \eqref{nontr1} are equivalent to those in \eqref{nontr} while the generalized Euler equation \eqref{eul} reduces to the Lagrangian inclusion \eqref{lagr}.

To verify further the complementary slackness conditions in \eqref{comp}, fix $i\in\{1,\ldots,m\}$ and suppose that $\ph_i(\ox)<0$. Then the continuity of $\ph_i$ at $\ox$ ensures that the pair $(\ox,0)$ is an interior point of the epigraphical set $\epi\ph_i$. It readily implies that $N((\ox,0);\epi\ph_i)=(0,0)$, and hence $\lm_i=0$. This yields $\lm_i\ph_i(\ox)=0$, which justifies \eqref{comp}.

To complete the proof of the theorem, it remains to check that the validity of \eqref{mfcq} ensures that $\lm_0=1$ in \eqref{lagr}. We easily arrive at this assertion while arguing by contradiction. $\h$\vspace*{0.05in}

If the constraint functions $\ph_i$, $i=1,\ldots,m$, are smooth around the reference point $\ox$, condition \eqref{mfcq} clearly reduces to the classical Mangasarian-Fromovitz constraint qualification. It suggests us to label this condition \eqref{mfcq} as the {\em generalized Mangasarian-Fromovitz constraint qualification}, or the {\em generalized MFCQ}.

\section{Variational Approach to Bilevel Optimization}\label{sec:bilevel}
\setcounter{equation}{0}

This section is devoted to describing some models of {\em bilevel programming} and a {\em variational approach} to them that involves {\em  nondifferentiable} optimal value functions in single-level problems of parametric optimization.

Let us first consider the following problem of {\em parametric optimization} with respect to the decision variable $y\in\R^m$ under each fixed parameter $x\in\R^n$:
\begin{equation}\label{lower}
\mbox{minimize }\ph(x,y)\;\mbox{ subject to }\;y\in F(x)\;\mbox{ with fixed }\;x\in\R^n,
\end{equation}
where $\ph\colon\R^n\times\R^m\to\R$ is the cost function and $F\colon\R^n\to\R^m$ is the constraint mapping in \eqref{lower}, which is called the {\em lower-level problem} of parametric optimization. Denoting by
\begin{equation}\label{sol-map}
S(x):=\mbox{argmin}\big\{\ph(x,y)\big|\;y\in F(x)\big\}
\end{equation}
the parameterized {\em solution set} for \eqref{lower} for each $x\in\R^n$ and given yet another cost function $\psi\colon\R^n\to\R^m$, we consider the {\em upper-level} parametric optimization problem of minimizing $\psi(x,y)$ over the lower-level solution map $S\colon\R^n\tto\R^m$ from \eqref{sol-map} written as:
\begin{equation}\label{upper}
\mbox{minimize }\;\psi(x,y)\;\mbox{ subject to }\;y\in S(x)\;\mbox{ for each }\;x\in\R^n.
\end{equation}
The {\em optimistic bilevel programming} model is defined by
\begin{equation}\label{opt}
\mbox{minimize }\;\mu(x)\;\mbox{ subject to }\;x\in\O,\;\mbox{ where }\;\mu(x):=\inf\big\{\psi(x,y)\big|\;y\in S(x)\big\},
\end{equation}
and where $\O\subset\R^n$ be a given constraint set. On the other hand, the {\em pessimistic bilevel programming} model is defined as follows:
\begin{equation}\label{pes}
\mbox{minimize }\;\eta(x)\;\mbox{ subject to }\;x\in\O,\;\mbox{ where }\;\eta(x):=\sup\big\{\psi(x,y)\big|\;y\in S(x)\big\}.
\end{equation}
We refer the reader to \cite{cms07}, \cite{demp03}--\cite{dmz19}, \cite{hen-sur11}, \cite{m18}, \cite{ye-zhu10}, \cite{zem16}, \cite{zhang03}, and the bibliographies therein for more details on both optimistic and pessimistic versions in bilevel programming, their local and global solutions as well as reformulations, modifications and relationships with other classes of optimization problems, theoretical and numerical developments, and various applications in finite-dimensional spaces. Investigations and applications of bilevel optimization problems in infinite dimensions can be found, e.g., in \cite{bgm07,benita-dem-meh16,dmn10,lig-morgan17,mnp11,zas12}.

The main attention in this and subsequent sections is paid to the application of the machinery and results of variational analysis and generalized differentiation to problems of bilevel optimization by implementing the {\em value function approach}. This approach is based on reducing bilevel programs to single-level problems of mathematical programming by using the nonsmooth {\em optimal value function} $\vt(x)$ of the lower-level problem defined in \eqref{marg-def}. Such a device was initiated by Outrata \cite{out90} for a particular class of bilevel optimization problems and was used by him for developing a numerical algorithm to solve bilevel programs. Then this approach was strongly developed by Ye and Zhu \cite{ye-zhu95} who employed it to derive necessary optimality conditions for optimistic bilevel programs by using Clarke's generalized gradients of optimal valued functions. More advanced necessary optimality conditions for optimistic bilevel programs in terms of the author's generalized differentiation reviewed in Section~\ref{basic} were developed in \cite{demp-dut-mor07,demp-mor-zem12,m18,mnp11,zem16}. Optimality and stability conditions for pessimistic bilevel models were derived in \cite{demp-mor-zem14,dmz19}.\vspace*{0.05in}

We restrict ourselves in what follows to implementing variational analysis and the aforementioned machinery of generalized differentiation within the value function approach to optimistic bilevel models with Lipschitzian data in finite-dimensional spaces. This allows us to most clearly communicate the basic variational ideas behind this approach, without additional technical complications. The variational results presented in the previous sections make our presentation self-contained and complete.

For simplicity we consider the optimistic bilevel model \eqref{opt} with only {\em inequality constraints} on the lower and upper levels described by
\begin{equation}\label{f}
F(x):=\big\{y\in\R^m\big|\;f_i(x,y)\le 0\;\mbox{ for }\;i=1,\ldots,r\big\},
\end{equation}
\begin{equation}\label{g}
\O:=\big\{x\in\R^n\big|\;g_j(x)\le 0\;\mbox{ for }\;j=1,\ldots,s\big\}.
\end{equation}
The reduction of the bilevel program \eqref{opt} with the constraints \eqref{f} and \eqref{g} to a single-level problem of nondifferentiable programming and deriving in this way necessary optimality conditions for it are given in the next section.

\section{Optimality Conditions for Lipschitzian Bilevel Programs}\label{optimality}
\setcounter{equation}{0}

The optimal value function \eqref{marg-def} for the lower-level program \eqref{lower} with the inequality constraints specified in \eqref{f} reads as
\begin{equation}\label{value}
\vt(x)=\inf\big\{\ph(x,y)\big|\;f_i(x,y)\le 0,\;i=1,\ldots,r\big\},\quad x\in\R^n.
\end{equation}
With the upper-level cost function $\psi$ given in \eqref{upper} and the upper-level constraints taken from \eqref{g}, consider the following single-level mathematical program with inequality constraints:
\begin{equation}\label{bilevel}
\begin{array}{ll}
&\mbox{minimize }\;\psi(x,y)\;\mbox{ subject to }\;g_j(x)\le 0,\;j=1,\ldots,s,\\
&f_i(x,y)\le 0,\;i=1,\ldots,r,\;\mbox{ and }\;\ph(x,y)\le\vt(x).
\end{array}
\end{equation}
We can easily observe that global optimal solutions to \eqref{bilevel} agree with those to problem \eqref{opt}, \eqref{f}, and \eqref{g}. Although this is not always the case for local minimizers, it is not hard to check that the local solutions to these optimization problems are also the same under the inner semicontinuity assumption on the solution map \eqref{sol-map} imposed in both main theorems obtained in this and next sections. To deriving further necessary optimality conditions for optimistic bilevel programming, we can therefore concentrate on the single-level optimization problem \eqref{bilevel}.

Looking at \eqref{bilevel}, observe that this problem is of type \eqref{mp} for which necessary optimality conditions are given in Theorem~\ref{lip-mp}, provided that all the functions involved are locally Lipschitzian. However, the direct application of Theorem~\ref{lip-mp} to problem \eqref{bilevel} is not efficient due to the structure of the last constraint therein defined via the lower-level optimal value function \eqref{value}. Indeed, it has been realized in bilevel programming that this constraint prevents the fulfillment of conventional constraint qualifications; in particular, the generalized MFCQ \eqref{mfcq}. To avoid this obstacle, Ye and Zhu \cite{ye-zhu95} introduced the following property postulating an appropriate behavior of the cost function in \eqref{bilevel} with respect to linear perturbations of the constraint $\ph(x,y)\le\vt(x)$. Consider the problem:
\begin{equation*}
\begin{array}{ll}
&\mbox{minimize }\;\psi(x,y)\;\mbox{ subject to }\;g_j(x)\le 0,\;j=1,\ldots,s,\\
&f_i(x,y)\le 0,\;i=1,\ldots,r,\;\mbox{ and }\;\ph(x,y)-\vt(x)+\nu=0\;\mbox{ as }\;\nu\in\R.
\end{array}
\end{equation*}

\begin{Def}\label{pcalm} Problem \eqref{bilevel} is {\sc partially calm} at its feasible solution $(\ox,\oy)$ if there exist a constant $\kappa>0$ and a neighborhood $U$ of $(\ox,\oy,0)\in\R^n\times\R^m\times\R$ such that
\begin{eqnarray}\label{calm}
\psi(x,y)-\psi(\ox,\oy)+\kappa|\nu|\ge 0
\end{eqnarray}
for all the triples $(x,y,\nu)\in U$ feasible to \eqref{bilevel}.
\end{Def}

There are various efficient conditions, which ensure the fulfillment of the partial calmness property for \eqref{bilevel}. They include the uniform sharp minimum condition \cite{ye-zhu95}, linearity of the lower-level problem with respect to the decision variable \cite{demp-zem13}, the kernel condition \cite{m18}, etc. On the other hand, partial calmness may fail in rather common situations; see \cite{m18} for more results and discussions on partial calmness and related properties.

The main impact of partial calmness to deriving necessary optimality conditions for \eqref{bilevel} is its equivalence to the possibility of transferring the troublesome constraint $\ph(x,y)\le\vt(x)$ into the {\em penalized} cost function as in the following proposition.

\begin{Prop}\label{calm-penal} Let $(\ox,\oy)$ be a partially calm feasible solution to problem \eqref{bilevel} with $\psi$ being continuous at this point. Then $(\ox,\oy)$ is a local optimal solution to the penalized problem
\begin{eqnarray}\label{pen}
\begin{array}{ll}
\mbox{minimize }\;\psi(x,y)+\kappa\big(\ph(x,y)-\vt(x)\big)\;\mbox{ subject to }\\
g_j(x)\le 0,\;j=1,\ldots,s,\;\mbox{ and }\;f_i(x,y)\le 0,\;i=1,\ldots,r,
\end{array}
\end{eqnarray}
where $\kappa>0$ is taken from \eqref{calm}. Conversely, any local optimal solution $(\ox,\oy)$ to \eqref{pen} with some number $\kappa>0$ is partially calm in \eqref{bilevel}.
\end{Prop}
{\bf Proof.} Taking $\kappa$ and $U$ from Definition~\ref{pcalm} and using the continuity of $\psi$ at $(\ox,\oy)$, we find $\gg>0$ and $\eta>0$ with $\Tilde U:=[(\ox,\oy)+\eta\B]\times(-\gg,\gg)\subset U$ and
\begin{eqnarray*}
|\psi(x,y)-\psi(\ox,\oy)|\le\kappa\gg\;\mbox{ for all }\;(x,y)-(\ox,\oy)\in\eta\B.
\end{eqnarray*}
Let us employ it to verifying that
\begin{eqnarray}\label{calm1}
\psi(x,y)-\psi(\ox,\oy)+\kappa\big(\ph(x,y)-\vt(x)\big)\ge 0\;\mbox{ whenever }\;g_j(x)\le 0,\;j=1,\ldots,s,
\end{eqnarray}
and $(x,y)\in[(\ox,\oy)+\eta\B]\cap\gph F$ with $F$ taken from \eqref{f}. If $(x,y,\vt(x)-\ph(x,y))\in\Tilde U$, then \eqref{calm1} follows from \eqref{calm}. In the remaining case where $(x,y,\vt(x)-\ph(x,y))\notin\Tilde U$, we get that
\begin{equation*}
\ph(x,y)-\vt(x)\ge\gg,\;\mbox{ and hence }\;\kappa\big(\ph(x,y)-\vt(x)\big)\ge\kappa\gg,
\end{equation*}
which also yields \eqref{calm1} by $\psi(x,y)-\psi(\ox,\oy)\ge-\kappa\gg$. The feasibility of $(\ox,\oy)$ to \eqref{bilevel} tells us that $\ph(\ox,\oy)-\vt(\ox)=0$, which verifies the first statement of the proposition. Arguing by contraction, we deduce the converse assertion directly from the definitions. $\h$\vspace*{0.05in}

Thus the imposed partial calmness allows us to deduce the original problem of optimistic bilevel optimization to the single-level mathematical program \eqref{pen} with conventional inequality constraints, where the troublesome term $\ph(x,y)-\vt(x)$ enters the penalized cost function. To derive necessary optimality conditions for \eqref{pen}, let us reformulate the generalized MFCQ \eqref{mfcq} in conventional bilevel terms as in the case of bilevel programs with smooth data \cite{demp03}.

We say that $(\ox,\oy)\in\R^n\times\R^m$ is {\em lower-level regular} if it satisfies the generalized MFCQ in the lower-level problem \eqref{lower}. This means due to the structure of \eqref{lower} that
\begin{equation}\label{lower-reg}
\Big[\disp\sum_{i\in I(\ox,\oy)}\lm_i v_i=0,\;\lm_i\ge 0\Big]\Longrightarrow\Big[\lm_i=0\;\mbox{ for all }\;i\in I(\ox,\oy)\Big]
\end{equation}
whenever $(u_i,v_i)\in\partial f_i(\ox,\oy)$ with some $u\in\in\R^n$ and $I(\ox,\oy):=\big\{i\in\{1,\ldots,r\}\big|\;f_i(\ox,\oy)=0\big\}$. Similarly, a point $\ox\in\R^n$ satisfying the upper-level constraints in \eqref{g} is {\em upper-level regular} if
\begin{equation}\label{upper-reg}
\Big[\disp 0\in\sum_{j\in J(\ox)}\lm_j\partial g_j(\ox),\;\lm_j\ge 0\Big]\Longrightarrow\Big[\lm_j=0\;\mbox{ whenever }\;j\in J(\ox)\Big]
\end{equation}
with the active constraint indexes $J(\ox):=\big\{j\in\{1,\ldots,s\}\big|\;g_j(\ox)=0\big\}$.\vspace*{0.05in}

Now we are ready for the application of Theorem~\ref{lip-mp} to the optimistic bilevel program in the equivalent form \eqref{pen}. To proceed, we have to verify first that the optimal value function $\vt(x)$ in the cost function of \eqref{pen} is locally Lipschitzian around the reference point and then to be able to upper estimate the basic subdifferential \eqref{bs} of the function $-\vt(\cdot)$. Since the basic subdifferential $\partial\vt(\ox)$ does not possess the plus-minus symmetry while its convex hull ${\rm co}\,\partial\vt(\ox)$ does, we need to convexify in the proof the set on the right-hand side of the upper estimate in \eqref{marg-bs}. In this way we arrive at the following major result.

\begin{Theorem}\label{bilevel-nc1} Let $(\ox,\oy)$ be a local optimal solution to the optimistic bilevel program in the equivalent form \eqref{bilevel} with the lower-level optimal value function $\vt(x)$ defined in \eqref{value}. Assume that all the functions $\ph,\psi,f_i,g_j$ are locally Lipschitzian around the reference point, that the lower-level solution map $S$ from \eqref{sol-map} is inner semicontinuous at $(\ox,\oy)$, that the lower-level regularity \eqref{lower-reg} and upper-level regularity \eqref{upper-reg} conditions hold, and that problem \eqref{bilevel} is partially calm at $(\ox,\oy)$ with constant $\kappa>0$. Then there exist multipliers $\lm_1,\ldots,\lm_r$, $\mu_1,\ldots,\mu_s$, and $\nu_1,\ldots,\nu_r$ satisfying the sign and complementary slackness conditions
\begin{equation}\label{lm}
\lm_i\ge 0,\;\;\lm_i f_i(\ox,\oy)=0\;\mbox{ for }\;i=1,\ldots,r,
\end{equation}
\begin{equation}\label{mu}
\mu_j\ge 0,\;\;\mu_j g_j(\ox)=0\;\mbox{ for }\;j=1,\ldots,s
\end{equation}
\begin{equation}\label{nu}
\nu_i\ge 0,\;\;\nu_i f_i(\ox,\oy)=0\;\mbox{ for }\;i=1,\ldots,r,
\end{equation}
together with the following relationships, which involve some vector $u\in\co\partial\vt(\ox)$:
\begin{equation}\label{lagr1}
(u,0)\in\co\partial\ph(\ox,\oy)+\disp\sum_{i=1}^r\nu_i\co\partial f_i(\ox,\oy),
\end{equation}
\begin{equation}\label{lagr2}
(u,0)\in\partial\ph(\ox,\oy)+\kappa^{-1}\partial\psi(\ox,\oy)+\disp\sum_{i=1}^r\lm_i\partial f_i(\ox,\oy)+\sum_{j=1}^s\mu_j\Big(\partial g_j(\ox),0\Big).
\end{equation}
\end{Theorem}
{\bf Proof.} It follows from Proposition~\ref{calm-penal} that $(\ox,\oy)$ is a local minimizer of the mathematical program \eqref{pen} with inequality constraints. To show that it belongs to problems of Lipschitzian programming considered in Section~\ref{program}, we need to check that the optimal value function \eqref{value} is locally Lipschitzian around $\ox$ under the assumptions made. This function is clearly l.s.c.\ around $\ox$, and thus its Lipschitz continuity around this point is equivalent to the condition $\partial^\infty\vt(\ox)=\{0\}$. It follows from the upper estimate \eqref{marg-ss} of Theorem~\ref{marg-sub} and the assumed inner semicontinuity of the solution map \eqref{sol-map} at $(\ox,\oy)$ that
\begin{equation*}
\partial^\infty\vt(\ox)\subset D^*F(\ox)(0)\;\mbox{ with }\;F(x)=\big\{y\in\R^m\big|\;f_i(x,y)\le 0,\;i=1,\ldots,r\big\}.
\end{equation*}
Corollary~\ref{marg-lip1} tells us that the local Lipschitz continuity of $\vt$ around $\ox$ follows from the Lipschitz-like property of the mapping $F\colon\R^n\tto\R^m$. Since
\begin{equation*}
\gph F=\big\{(x,y)\in\R^n\times\R^m\big|\;f_i(x,y)\le 0,\;i=1,\ldots,r\big\},
\end{equation*}
we deduce that $D^*F(\ox,\oy)(0)=\{0\}$ from the assumed lower-level regularity due to the coderivative definition and the normal come intersection rule in Theorem~\ref{nor-inter}. Thus $F$ is Lipschitz-like around $(\ox,\oy)$ by the coderivative criterion \eqref{cc}, and the Lipschitz continuity of $\vt(\cdot)$ is verified.

Now we can apply to problem \eqref{pen} the necessary optimality conditions for Lipschitzian programs obtained in Theorem~\ref{lip-mp}. The assumed lower-level regularity and upper-level regularity clearly imply that the generalized MFCQ condition \eqref{mfcq} holds. Thus there are multiplies $\lm_1,\ldots,\lm_r$ and $\mu_1,\cdots,\mu_s$ satisfying the sign and complementary slackness conditions in \eqref{lm} and \eqref{mu} for which we have
\begin{equation}\label{lagr3}
\begin{array}{ll}
0&\in\partial\psi(\ox,\oy)+\kappa\partial\ph(\ox,\oy)+\Big(\kappa\partial(-\vt)(\ox),0\Big)\\
&+\disp\sum_{i=1}^r\lm_i\partial f_i(\ox,\oy)+\sum_{j=1}^s\mu_j\Big(\partial g_j(\ox),0\Big).
\end{array}
\end{equation}
To estimate $\partial(-\vt)(\ox)$ in \eqref{lagr3}, recall that
\begin{equation*}
\partial(-\vt)(\ox)\subset\Bar\partial(-\vt)(\ox)=-\Bar\partial\vt(\ox)=-\co\partial\vt(\ox),
\end{equation*}
where $\Bar\partial$ stands for Clarke's generalized gradient of locally Lipschitzian functions that possesses the plus-minus symmetry property \cite{cl}. Using it in \eqref{lagr3}, we get $u\in\co\partial\vt(\ox)$ such that
\begin{equation}\label{lagr4}
\kappa(u,0)\in\partial\psi(\ox,\oy)+\partial\ph(\ox,\oy)+\disp\sum_{i=1}^r\lm_i\partial f_i(\ox,\oy)+\sum_{j=1}^s\Big(\mu_j\partial g_j(\ox),0\Big).
\end{equation}
Applying the convexified subdifferential estimates \eqref{marg-bs} from Theorem~\ref{marg-sub} to the optimal value function \eqref{value} allows us to find multipliers $\nu_1,\ldots,\nu_r$ satisfying the sign and complementary slackness conditions in \eqref{nu} that ensure the validity of \eqref{lagr1}. To verify finally \eqref{lagr2}, we divide \eqref{lagr4} by $\kappa>0$ with keeping the same notation for the scaled multipliers $\lm_i$ and $\mu_j$. $\h$\vspace*{0.05in}

The next section presents an independent set of necessary optimality conditions for optimistic bilevel programs with Lipschitzian data that are obtained without using any convexification while employing instead yet another variational device and subdifferential calculus rule.

\section{Bilevel Optimization via Subdifferential Difference Rule}\label{difference}
\setcounter{equation}{0}

Considering the single-level problem \eqref{pen}, which we are finally dealing with while deriving necessary optimality conditions for optimistic bilevel programs, note that the objective therein contains the {\em difference} of two nonsmooth functions. The basic subdifferential \eqref{bs} does not possesses any special rule for difference of nonsmooth functions, but the regular subdifferential \eqref{rs} does, as was first observed in \cite{mor-nam-yen} by using a smooth variational description of regular subgradients. Here we employ this approach to establish necessary optimality conditions for Lipschitzian bilevel programs that are different from those in Theorem~\ref{bilevel-nc1}.\vspace*{0.05in}

The derivation of these necessary optimality conditions are based on the following two results, which are certainly of their independent interest. The first one provides a {\em smooth variational description} of regular subgradients of arbitrary functions $\ph\colon\R^n\to\oR$.

\begin{Lemma}\label{smooth-reg} Let $\ph\colon\R^n\to\oR$ be finite at $\ox$, and let $v\in\Hat\partial\ph(\ox)$. Then there exists a neighborhood $U$ of $\ox$ and a function $\psi\colon U\to\R$ such that $\psi(\ox)=\ph(\ox)$, that $\psi$ is Fr\'echet differentiable at $\ox$ with $\nabla\psi(\ox)=v$, and that the difference $\psi-\ph$ achieves at $\ox$ its local maximum on $U$.
\end{Lemma}
{\bf Proof.} We proceed geometrically due to the relationship
\begin{equation*}
\Hat\partial\ph(\ox)=\big\{v\in\R^n\big|\;(v,-1)\in\Hat N\big((\ox,\ph(\ox));\epi\ph\big)\big\},
\end{equation*}
which reduces the claimed assertion to the following: $v\in\Hat N(\ox;\O$ if and only if there exists a neighborhood $U$ of $\ox$ and a function $\psi\colon\R^n\to\R$ such that $\psi$ is Fr\'echet differentiable at $\ox$ with $\nabla\psi(\ox)=v$ while achieving at $\ox$ its local maximum relative to $\O$.

To verify the latter, observe that for any $\psi\colon U\to\R$ satisfying the listed properties we get
\begin{equation*}
\psi(x)=\psi(\ox)+\la v,x-\ox\ra+o(\|x-\ox\|)\le\psi(\ox)\;\mbox{ whenever }\;x\in U.
\end{equation*}
It shows that $\la v,x-\ox\ra+o(\|x-\ox\|)\le 0$, and thus $v\in\Hat N(\ox;\O)$ by definition \eqref{pre}. Conversely, pick $v\in\Hat N(\ox;\O)$ and define the function
\begin{eqnarray*}
\psi(x):=\left\{\begin{array}{ll}
\min\big\{0,\la v,x-\ox\ra\big\}&\mbox{ if }\;x\in\O,\\
\la v,x-\ox\ra&\mbox{ otherwise}.
\end{array}\right.
\end{eqnarray*}
It is easy to check that this function enjoys all the properties listed above. $\h$\vspace*{0.05in}

The second lemma gives us the aforementioned {\em difference rule} for regular subgradients.

\begin{Lemma}\label{diff-rule} Consider two arbitrary functions $\ph_1,\ph_2\colon\R^n\to\oR$ that are finite at $\ox$ and assume that $\Hat\partial\ph_2(\ox)\ne\emp$. Then we have the inclusions
\begin{eqnarray}\label{H1}
\Hat\partial(\ph_1-\ph_2)(\ox)\subset\bigcap_{v\in\Hat\partial\ph_2(\ox)}\Big[\Hat\partial\ph_1(\ox)-v\Big]\subset
\Hat\partial\ph_1(\ox)-\Hat\partial\ph_2(\ox).
\end{eqnarray}
It implies, in particular, that any local minimizer $\ox$ of the difference function $\ph_1-\ph_2$ satisfies the necessary optimality condition
\begin{eqnarray}\label{H2}
\Hat\partial\ph_2(\ox)\subset\Hat\partial\ph_1(\ox).
\end{eqnarray}
\end{Lemma}
{\bf Proof.} Starting with (\ref{H1}), pick $v\in\Hat\partial\ph_2(\ox)$. Then the smooth variational description of $v$ from Lemma~\ref{smooth-reg} gives us $\psi\colon U\to\R$ on a neighborhood $U$ of $\ox$ that is differentiable at $\ox$ with
\begin{eqnarray*}
\psi(\ox)=\ph_2(\ox),\quad\nabla\psi(\ox)=v,\;\mbox{ and }\;\psi(x)\le\ph_2(x)\;\mbox{ whenever }\;x\in U.
\end{eqnarray*}
Fix further an arbitrary vector $w\in\Hat\partial(\ph_1-\ph_2)(\ox)$ and for any $\ve>0$ find $\gg>0$ such that
\begin{eqnarray*}
\begin{array}{ll}
\la w,x-\ox\ra&\le\ph_1(x)-\ph_2(x)-\big(\ph_1(\ox)-\ph_2(\ox)\big)+\ve\|x-\ox\|\\
&\le\ph_1(x)-\psi(x)-\big(\ph_1(\ox)-\psi(\ox)\big)+\ve\|x-\ox\|
\end{array}
\end{eqnarray*}
if $\|x-\ox\|\le\gg$. Due to the differentiability of $\psi$ at $\ox$ we deduce that
\begin{equation*}
w\in\Hat\partial(\ph_1-\psi)(\ox)=\Hat\partial\ph_1(\ox)-\nabla\psi(\ox)=\Hat\partial\ph_1(\ox)-v
\end{equation*}
and thus verify both inclusions in \eqref{H1}.

Observing that \eqref{H2} is trivial if $\Hat\partial\ph_2(\ox)=\emp$, assume the opposite and fix any $v\in\Hat\partial\ph_2(\ox)$. Then it follows from \eqref{H1} and the obvious Fermat stationary rule via regular subgradients that
\begin{equation*}
0\in\Hat\partial(\ph_1-\ph_2)(\ox)\subset\Hat\partial\ph_1(\ox)-v.
\end{equation*}
It shows that $v\in\Hat\partial\ph_1(\ox)$ and thus verifies the fulfillment of \eqref{H2}. $\h$\vspace*{0.05in}

Now we are in a position to derive refined necessary optimality conditions for optimistic bilevel programs with Lipschitzian data.

\begin{Theorem}\label{bilevel-nc2} Let $(\ox,\oy)$ be a local optimal solution to the optimistic bilevel program in the equivalent form \eqref{bilevel} without upper level constraints. Suppose that all the functions $\ph,\psi,f_i$ are locally Lipschitzian around the reference point, that the lower-level solution map $S$ in \eqref{sol-map} is inner semicontinuous at $(\ox,\oy)$, that the lower-level regularity \eqref{lower-reg} condition holds, and that problem \eqref{bilevel} is partially calm at $(\ox,\oy)$ with constant $\kappa>0$. Assume in addition that $\Hat\partial\vt(\ox)\ne\emp$ for the optimal value function \eqref{value}. Then there exists a vector $u\in\Hat\partial\vt(\ox)$ together with multipliers $\lm_1,\ldots,\lm_r$ and $\nu_1,\ldots,\nu_r$ for $i=1,\ldots,r$ satisfying the sign and complementary slackness conditions in \eqref{nu} and \eqref{lm}, respectively, such that
\begin{equation}\label{lagr11}
(u,0)\in\partial\ph(\ox,\oy)+\disp\sum_{i=1}^r\nu_i\partial f_i(\ox,\oy),
\end{equation}
\begin{equation}\label{lagr21}
(u,0)\in\partial\ph(\ox,\oy)+\kappa^{-1}\partial\psi(\ox,\oy)+\disp\sum_{i=1}^r\lm_i\partial f_i(\ox,\oy).
\end{equation}
\end{Theorem}
{\bf Proof.} We get from the partial calmness penalization in Proposition~\ref{calm-penal} employed together with the infinite penalization of the lower-level constraints that $(\ox,\oy)$ a local minimizer for the unconstrained optimization problem
\begin{equation}\label{lip1}
\mbox{minimize }\;\psi(x,y)+\kappa\big(\ph(x,y)-\vt(x)\big)+\delta\big((x,y);\gph F\big)
\end{equation}
with the mapping $F\colon\R^n\tto\R^m$ defined in \eqref{f}. Then applying to \eqref{lip1} the difference rule \eqref{H2} from Proposition~\ref{diff-rule} gives us the inclusion
\begin{equation}\label{lip2}
\big(\kappa\Hat\partial\vt(\ox),0\big)\subset\Hat\partial\big(\psi(\cdot)+\kappa\ph(\cdot)+\delta(\cdot;\gph F)\big)(\ox,\oy).
\end{equation}
At the same time it follows from the proof of Theorem~\ref{marg-sub} that
\begin{equation}\label{lip3}
\big(\Hat\partial\vt(\ox),0\big)\subset\Hat\partial\big(\ph(\cdot)+\delta\big(\cdot;\gph F\big)\big)(\ox,\oy).
\end{equation}
Replacing $\Hat\partial$ by the larger $\partial$ on the right-hand sides of \eqref{lip2} and \eqref{lip3} and then using the basic subdifferential sum rule from Theorem~\ref{sub-sum} yield the inclusions
\begin{equation}\label{lip4}
\begin{aligned}
\Big(\kappa\Hat\partial\vt(\ox),0\Big)&\subset
\partial\psi(\ox,\oy)+\kappa\partial\ph(\ox,\oy)+N\big((\ox,\oy);\gph F\big),\\
\Big(\Hat\partial\vt(\ox),0\Big)&\subset\partial\ph(\ox,\oy)+N\big((\ox,\oy);\gph F\big).
\end{aligned}
\end{equation}
Employing in these inclusions Theorem~\ref{nor-inter} under the lower-level regularity of $(\ox,\oy)$ with the usage of the singular subdifferential characterization of Lipschitzian functions in \eqref{lip-func}, we get
\begin{equation*}
N\big((\ox,\oy);\gph F\big)\subset\bigcup\Big\{\disp\sum_{i=1}^r\lm_i\partial f_i(\ox,\oy)\Big|\;\lm_i\ge 0,\;\lm_i f_i(\ox,\oy)=0\;\mbox{ as }\; i=1,\ldots,r\Big\}.
\end{equation*}
It allows us to deduce from \eqref{lip4} the existence of $u\in\Hat\partial\vt(\ox)$ ensuring the validity of \eqref{lagr11} and
\begin{equation*}
\kappa(u,0)\in\partial\psi(\ox,\oy)+\kappa\partial\ph(\ox,\oy)+\disp\sum_{i=1}^r\lm_i\partial f_i(\ox,\oy).
\end{equation*}
Dividing the latter by $\kappa>0$, we arrive at \eqref{lagr21} and thus complete the proof of the theorem. $\h$\vspace*{0.05in}

Observe that we always have $\Hat\partial\vt(\ox)\ne\emp$ if the optimal value function \eqref{value} is {\em convex}, which is surely the case when all the functions $\ph$ and $f_i$ therein are convex. If in addition the upper-level data are also convex, more general results were derived in \cite{dmn10} for problems of {\em semi-infinite programming} with arbitrarily number of inequality constraints in locally convex topological vector spaces by reducing them to problems of {\em DC programming} with objectives represented as {\em differences of convex} functions. Note further that the necessary optimality conditions obtained in Theorems~\ref{bilevel-nc1} and \ref{bilevel-nc2} are independent of each other even in the case of bilevel programs with smooth data. In particular, we refer the reader to \cite[Example~6.24]{m18} for illustrating this statement and for using the obtained results to solve smooth bilevel programs. Finally, we mention the possibility to replace the inner semicontinuity assumption on the solution map $S(x)$ imposed in both Theorems~\ref{bilevel-nc1} and \ref{bilevel-nc2} by the {\em uniform boundedness} of this map in finite-dimensions, or by its {\em inner semicompactness} counterpart in infinite-dimensional spaces; cf.\ \cite{demp-dut-mor07,dmz19,m18,mnp11} for similar transitions in various bilevel settings.

\section{Concluding Remarks and Open Questions}\label{conclusion}
\setcounter{equation}{0}

In this self-contained chapter of the book we described a variational approach to bilevel optimization with its implementation to deriving advanced necessary optimality conditions for optimistic bilevel programs in finite-dimensional spaces. The entire machinery of variational analysis and generalized differentiation (including the fundamental extremal principle, major calculus rules, and subdifferentiation of optimal value functions), which is needed for this device, is presented here with the proofs. The given variational approach definitely has strong perspectives for further developments. Let us briefly discuss some open questions in this direction.

$\bullet$ The major difference between the optimistic model \eqref{opt} and pessimistic model \eqref{pes} in bilevel programming is that the latter invokes the {\em supremum} marginal/optimal value function instead of the infimum type in \eqref{opt}. Subdifferentiation of the supremum marginal functions is more involved in comparison with that of the infimum type. Some results in this vein for problems with Lipschitzian data can be distilled from the recent papers \cite{mor-nghia13,mor-nghia14,perez}, while their implementation in the framework of pessimistic bilevel programs is a challenging issue.

$\bullet$ The given proof of the necessary optimality conditions for bilevel programs in Theorem~\ref{bilevel-nc1} requires an upper estimate of $\partial(-\vt)(\ox)$, which cannot be directly derived from that for $\partial\vt(\ox)$ since $\partial(-\vt)(\ox)\ne-\partial\vt(\ox)$. To obtain such an estimate, we used the subdifferential convexification and the fact that the convexified/Clarke subdifferential of Lipschitz continuous functions possesses the plus-minus symmetry. However, there is a nonconvex subgradient set that is much smaller than Clarke's one while having this symmetry. It is the {\em symmetric subdifferential} $\partial^0\vt(\ox):=\partial\vt(\ox)\cup(-\partial(-\vt)(\ox))$, which enjoys full calculus induced by the basic one. Efficient evaluations of $\partial^0\vt(\ox)$ for infimum and supremum marginal functions would lead us to refined optimality conditions for both optimistic and pessimistic models in bilevel optimization.

$\bullet$ The partial calmness property used in both Theorems~\ref{bilevel-nc1} and \ref{bilevel-nc2} seems to be rather restrictive when the lower-level problem is nonlinear with respect to the decision variable. It is a challenging research topic to relax this assumption and to investigate more the uniform weak sharp minimum property and its modifications that yield partial calmness.

$\bullet$  One of the possible ways to avoid partial calmness in bilevel programming is as follows. Having the solution map $S\colon\R^n\tto\R^m$ to the lower-level problem, consider the constrained upper-level problem given by
\begin{equation}\label{multi}
\mbox{minimize }\;\Psi(x):=\psi\big(x,S(x)\big)\;\mbox{ subject to }\;x\in\O,
\end{equation}
where $\psi\colon\R^n\times\R^m\to\oR$ is the cost function on the upper level with the upper-level constraint set $\O\subset\R^n$, and where the minimization of $\Psi\colon\R^n\tto\R$ is understood with respect to the standard order on $\R$. Then \eqref{multi} is a problem of {\em set-valued optimization} for which various necessary optimality conditions of the coderivative and subdifferential types can be found in \cite{m18} and the references therein. Evaluating the coderivatives and subdifferentials of the composition $\psi(x,S(x))$ in terms of the given lower-level and upper-level data of bilevel programs would lead us to necessary optimality conditions in both optimistic and pessimistic models. There are many open questions arising in efficient realizations of this approach for particular classes of problems in bilevel optimization even with smooth initial data. We refer the reader to the paper by Zemkoho \cite{zem16} for some recent results and implementations in this direction.

Note that a somewhat related approach to bilevel optimization was developed in \cite{bgm07}, where a lower-level problem was replaced by the corresponding KKT system described by a certain {\em generalized equation} of the Robinson type \cite{rob}. Applying to the latter necessary optimality conditions for upper-level problems with such constraints allowed us to establish verifiable results for the original nonsmooth problem of bilevel programming.

$\bullet$ Henrion and Surowiec suggested in \cite{hen-sur11} a novel approach to derive necessary optimality conditions for optimistic bilevel programs with ${\cal C}^2$-smooth data and convex lower-level problems. Their approach used a reduction to {\em mathematical programs with equilibrium constraints} (MPECs) and allowed them to significantly relax the partial calmness assumption. Furthermore, in this way they obtained new necessary optimality conditions for the bilevel programs under consideration, which are described via the Hessian matrices of the program data.

A challenging direction of the future research is to develop the approach and results from \cite{hen-sur11} to nonconvex bilevel programs with nonsmooth data. It would be natural to replace in this way the classical Hessian in necessary optimality conditions by the generalized one (known as the {\em second-order subdifferential}) introduced by the author in \cite{m92} and then broadly employed in variational analysis and its applications; see, e.g., \cite{m18} and the references therein.

$\bullet$ As has been long time realized, problems of bilevel optimization are generally {\em ill-posed}, which creates serious computational difficulties for their numerical solving; see, e.g., \cite{cms07,demp03,zem16} for more details and discussions. Furthermore, various regularization methods and approximation procedures devised in order to avoid ill-posedness have their serious drawbacks and often end up with approximate solutions, which may be far enough from optimal ones. Thus it seems appealing to deal with ill-posed bilevel programs how they are and to develop {\em numerical algorithms} based on the obtained necessary optimality conditions. Some results in this vein are presented in \cite{zem16} with involving necessary optimality conditions of the type discussed here, while much more work is required to be done in this very important direction with practical applications.

\end{document}